\documentclass[12pt]{article}
\usepackage{latexsym,amsfonts,amsthm,amsmath,amscd,amssymb}
\usepackage[dvips]{graphicx}

\newtheorem{lemma}{Lemma}[section]

\newtheorem{prop}[lemma]{Proposition}

\newtheorem{ques}[lemma]{Question}

\newcommand\matN{{\mathbb{N}}}

\newcommand\matP{{\mathbb{P}}}
\newcommand\matR{{\mathbb{R}}}

\newcommand\matZ{{\mathbb{Z}}}

\newcommand\hatM{{\widehat{M}}}
\newcommand\Mbar{{\overline{M}}}

\renewcommand{\hbar}{{\overline{h}}}

\newfont{\Got}{eufm10 scaled 1200}

\newcommand{\mycap} [1] {\caption{\footnotesize{#1}}}

\begin{document}

\title{Combinatorial and geometric\\ methods in topology}

\author{Carlo Petronio\\ with an appendix by\\ Damian Heard and Ekaterina
Pervova}

\maketitle

\begin{abstract}
\noindent Starting from the (apparently) elementary problem of
deciding how many different topological spaces can be obtained by
gluing together in pairs the faces of an octahedron, we will
describe the central r\^ole played by hyperbolic geometry within
three-dimensional topology. We will also point out the striking
difference with the two-dimensional case, and we will review some
of the results of the combinatorial and computational approach to
three-manifolds developed by different mathematicians over the
last several years.

\end{abstract}

\noindent The octahedron, denoted henceforth by $O$, is one of the
favourite toys of every geometer, being one of the five Platonic
solids. In this note we will investigate the following:

\begin{ques}\label{basic:ques}
How many different topological spaces can be obtained
by gluing together in
pairs the faces of $O$?
\end{ques}

Question~\ref{basic:ques} has a rather transparent combinatorial
flavour and appears to be well-suited to computer investigation,
but the complete answer would be extremely difficult to obtain
without the aid of some rather sophisticated geometric tools
developed over the last three decades by a number of
mathematicians. It is indeed mostly thanks to hyperbolic geometry
that one is able to show that certain gluings of $O$, despite
being very similar to each other under many respects, are in fact
distinct.

There is a very natural family of lesser siblings of
Question~\ref{basic:ques}, involving (two-dimensional) polygons
rather than a (three-dimensional) polyhedron (as the octahedron
is), and we will show below that the answers to these
lower-dimensional questions are sharply different from both a
qualitative and quantitative viewpoint. More precisely, it turns
out that identifying the spaces obtained by gluing together the
edges of a polygon is very easy, and that the number of possible
different results is very small if compared to the number of
(combinatorially inequivalent) gluing patterns. On the other hand,
in dimension three identifying the results is often only possible
using hyperbolic geometry, and there is a rather large variety of
different results. This can be viewed as a manifestation of the
crucial r\^ole played by hyperbolic geometry in the context of
three-dimensional topology, as chiefly witnessed by Thurston's
geometrization.

As opposed to looking at lower-dimensional analogues of
Question~\ref{basic:ques}, one can also view it as a special instance of a more
general family of three-dimensional problems,
where one considers larger polyhedra (or finite families
of polyhedra). These problems have attracted a considerable attention during the
last several years, and we will include a brief survey of the main
results obtained.

\vspace{1cm}

\noindent\textsc{Acknowledgements}. Part of this work was carried
out while the author was visiting the Universit\'e Paul Sabatier in
Toulose and the Columbia University in New York. He is grateful to
both these institutions for financial support, and he would like to
thank Michel Boileau and Dylan Thurston for their warm hospitality
and inspiring mathematical discussions. He is
also glad that Damian Heard and Katya Pervova
accepted to contribute to the paper with an appendix where they explain
how they have determined some of the
entries in Table~\ref{final:tab}.

\section{One dimension down}
Let us denote by $P_k$ a polygon with $2k$ edges. Since we are only interested in topology,
this makes sense for $k\geqslant 1$. In other words, $P_k$ is just the 2-dimensional closed
disc with boundary circle subdivided into $2k$ segments.

\begin{ques}\label{2dim:ques}
How many different topological spaces can be obtained
by gluing together in pairs the edges of $P_k$?
\end{ques}

We will investigate this question only with a restriction
on the gluings. Namely, we endow $P_k$ with some orientation, we give
its edges the induced orientation, and we require the gluing maps
to reverse this orientation of the edges.
With this proviso a gluing pattern is determined
by the instructions of which edges should get glued to which.
If one gives labels from $1$ to $2k$ to the edges of $P_k$
(as one should do if one wants to feed $P_k$ to a computer)
one easily sees that the number of different patterns is $(2k-1)!!$
This is because one must first choose which edge gets glued to edge number $1$,
and there are $2k-1$ choices for this. Then one must locate the first
edge which is left free after the first gluing, and select one
of the remaining $2k-3$ edges to glue it to, and so on.

As just mentioned, counting the different gluing patterns is easy,
but one must also note that $P_k$ has some symmetries, given by the
action of the dihedral group $D_{2k}$, which has $4k$ elements. So
the number of combinatorially inequivalent patterns is smaller than
$(2k-1)!!$, but at least $(2k-1)!!/4k$ (and actually always strictly
larger than this value, because some patterns will have symmetries
themselves). Giving a general formula for the number of inequivalent
patterns for arbitrary $k$ is probably rather complicated, but it is
not too difficult to write a piece of computer code that does the
counting for small $k$. The author has implemented this counting
algorithm using the language Haskell and computed the number of
gluings for $k\leqslant 6$, with results as illustrated below.

Knowing how many gluing patterns exist for some $k$, however, does
not answer Question~\ref{2dim:ques}, because inequivalent patterns
may well give homeomorphic spaces after gluing, and as a matter of
fact they very often do. But one can notice that the topological
space resulting from a gluing is always a surface, and an
orientable one under our restriction that the gluings should
reverse the orientation of the edges. Now, an orientable surface
$\Sigma$ is classified by its genus $g(\Sigma)$, an integer that
can attain any non-negative value, and the genus can easily be
computed from the Euler characteristic thanks to the relation
$\chi(\Sigma)=2(1-g(\Sigma))$. Computing $\chi(\Sigma)$ for a
surface $\Sigma$ obtained by gluing together the edges of $P_k$ is
also very easy, because $\chi(\Sigma)=v-k+1$, where $v$ is the
number of equivalence classes of vertices under the relation
generated by the gluings. This implies that $\chi(\Sigma)\geqslant
2-k$, so $g(\Sigma)\leqslant k/2$. On the other hand it is well
known that the surface of genus $g>0$ can be realized by gluing the
edges of $P_{2g}$, and as a matter of fact also all lower-genus
surfaces can.  This discussion implies the following:

\begin{prop}\label{2d:prop}
The number of combinatorially inequivalent patterns for the
edges of $P_k$ is at least $(2k-1)!!/4k$, while the number of
distinct resulting surfaces is $[k/2]+1$.
\end{prop}

\begin{table}
\begin{center}
\begin{tabular}{c||c|c|c|c|c|c}
$k$ & 1 & 2 & 3 & 4 & 5 & 6 \\ \hline\hline
$\#$(gluings) & 1 & 3 & 15 & 105 & 945 & 10395 \\ \hline
$\#$(inequivalent gluings)& 1 & 2 & 5 & 17 & 79 & 554 \\ \hline
$\#$(resulting surfaces) & 1 & 2 & 2 & 3 & 3 & 4 \\
\end{tabular}
\end{center}
\mycap{Numbers of gluing patterns for the $2k$-edged
polygon and of the resulting surfaces\label{2d:tab}}
\end{table}

The exact figures for $k\leqslant 6$ are given in
Table~\ref{2d:tab}. The qualitative conclusion we can draw from the
above discussion, and in particular from Proposition~\ref{2d:prop}
and Table~\ref{2d:tab}, is that \emph{in two dimensions it is easy
to recognize the surface resulting from a gluing of the edges of a
polygon, and the number of surfaces that arise is very small
compared to the number of inequivalent gluings}. We will see that
the situation is sharply different in
three dimensions.

\section{Three-manifolds (with boundary)}
In the previous section we have taken for granted that the result of
a gluing of the edges of a polygon is a surface, namely a (compact)
topological space locally homeomorphic to the plane $\matR^2$. Since
we will see that a similar result almost but not quite holds in
three dimensions, we will now
prove the 2-dimensional statement, in order to show where the
difficulty can arise. As we do this we momentarily drop the
restriction that the gluing should reverse the orientation of the
edges, because the conclusion holds in general.

\begin{prop}\label{surf:prop}
The space resulting from an arbitrary simplicial pairing
of the edges of the polygon $P_k$ is a surface.
\end{prop}

\begin{proof}
Denote by $X$ the result of the gluing and by $\pi:P_k\to X$ the projection.
There are three different sorts of points $x$ of $X$ that we must look at,
depending on what $\pi^{-1}(x)$ is. The possibilities are as follows:
\begin{itemize}
\item $\pi^{-1}(x)$ is a single interior point $p$ of $P_k$;
\item $\pi^{-1}(x)$ consists of two points $p_1$ and $p_2$ lying in the
interior of two edges of $P_k$;
\item $\pi^{-1}(x)$ consists of some vertices $v_1,\ldots,v_n$ of $P_k$.
\end{itemize}
For the first type, since $\pi$ is a local homeomorphism near $p$, it is
clear that $x$ is a surface point of $X$. For the second type, a neighbourhood
of $x$ in $X$ is obtained by gluing neighbourhoods of $p_1$ and $p_2$ in $P_k$.
These neighbourhoods can be chosen to be half-discs and the gluing identifies
their diameters, so again $x$ is a surface point. Getting to the last type
of points, we note that a neighbourhood of each $v_i$ in $P_k$
can be viewed as the cone with vertex $v_i$ over a segment.
Therefore a neighbourhood of $x$ in $X$ is the cone over the space obtained from the
union of $n$ segments by gluing together in pairs their endpoints. This space
is connected, so it is a circle, whence a neighbourhood of $x$ in $X$ is a disc,
which proves that $x$ is again a surface point.
\end{proof}

Turning to three dimensions and our octahedron $O$, we introduce
again the restriction that the gluings of the faces should reverse
the orientation induced by some orientation of $O$. For instance, if
we label by $v_0,\ldots,v_5$ the vertices of $O$ so that the faces
are $(v_i,v_j,v_{j+1})$ for $i=0,5$ and $j=1,\ldots,4$, with the
convention that $j+1=1$ for $j=4$, then the maps
$$(v_0,v_1,v_2)\ \longrightarrow (v_0,v_3,v_2),\qquad (v_0,v_1,v_2)\ \longrightarrow (v_5,v_1,v_2)$$
define gluings that
are acceptable for us, while
$$(v_0,v_1,v_2)\ \longrightarrow (v_0,v_2,v_3),\qquad (v_0,v_1,v_2)\ \longrightarrow (v_5,v_2,v_1)$$
do not. We also
recall that a 3-manifold is a space locally homeomorphic to
$\matR^3$, and we denote by $O^{(0)}$ the set of vertices of $O$.

\begin{prop}\label{3mfld:prop}
Let $\varphi$ be a simplicial orientation-reversing pairwise gluing of the faces of
the octahedron $O$. Define
$$X(\varphi)=O/\varphi,\qquad V(\varphi)=O^{(0)}/\varphi,\qquad Y(\varphi)=X(\varphi)\setminus V(\varphi).$$
Then $Y(\varphi)$ is a $3$-manifold.
\end{prop}

\begin{proof}
Let $\pi:O\to X(\varphi)$ be the projection. We first prove the
following claim: \emph{the restriction of $\pi$ to the interior of
each edge of $O$ is injective}. To show this, we assume $\pi$
identifies together some edges $e_0,\ldots,e_{m-1}$ of $O$. More
precisely, we assume $e_i$ is contained in some faces $F_i$ and
$F'_i$ of $O$ and that there is a gluing from $F'_i$ to $F_{i+1}$
mapping $e_i$ to $e_{i+1}$, with indices understood modulo $m$.
Furthermore we give to each $F_i$
and $F'_i$ the orientation induced by $O$, and to
each $e_i$ the orientation induced by $F_i$. Since $F'_i$ induces
on $e_i$ the negative orientation and the gluing from $F'_i$ to
$F_{i+1}$ reverses the orientation, we see that this gluing matches
the orientations of $e_i$ and $e_{i+1}$. Therefore the composition
of all the gluings $e_0\to e_1\to\ldots\to e_{m-1}\to e_0$ is an
orientation-preserving simplicial self-gluing of $e_0$. This gluing
is the identity, so $\pi$ restricted to the interior of $e_0$ is
injective. Our claim is proved.

As in the 2-dimensional case we now take $x\in Y(\varphi)$ and distinguish
according to what $\pi^{-1}(x)$ is. The possibilities are as follows:
\begin{itemize}
\item $\pi^{-1}(x)$ is a single interior point of $O$;
\item $\pi^{-1}(x)$ consists of two interior
points of some faces of $O$;
\item $\pi^{-1}(x)$ consists of some points
$p_1,\ldots,p_n$ lying in the interior of $n$ edges of $O$.
\end{itemize}
Again points of the first two types are obviously 3-manifold points.
For the last type we note that a neighbourhood of $p_i$ in $O$
can be viewed as the cone with vertex $p_i$ over the bigon $P_2$.
Therefore a neighbourhood of $x$ in $X(\varphi)$ is the cone over a connected
space obtained from a finite union of bigons by gluing together in
pairs their edges. Of course by performing all but one gluing we
can reduce to the case of a single bigon. Now, depending on whether
the gluing reverses or preserves the orientation induced on the edges
of the bigon, the resulting space is either the 2-sphere or a
non-orientable surface, called the projective plane. In the first case, since
the cone over a 2-sphere is a 3-disc, we see that $x$ is a 3-manifold point.
The same conclusion would be false in the second case, but an argument
very close to that presented above to establish the claim shows that the second
case actually never occurs under the current restriction that the
gluings should reverse the orientation. The proof is complete.
\end{proof}

We will now see that it is not possible to improve
Proposition~\ref{3mfld:prop} showing that the whole of $X(\varphi)$
is a 3-manifold. (We will later see that this is rather an advantage
than a disadvantage). To this end we define the \emph{truncated
octahedron} $O_{\rm t}$ as $O$ minus small open pyramidal
neighbourhoods of its vertices, chosen so that any simplicial gluing
of the faces of $O$ will match their triangular faces. Note that
$O_{\rm t}$ is a polyhedron bounded by 8 hexagons and 6 squares, and
that any simplicial gluing of the faces of $O$ induces a gluing
between the hexagonal faces of $O_{\rm t}$. Before stating our
result we recall that a 3-manifold \emph{with boundary} is a space
$M$ locally homeomorphic to a closed half-space of $\matR^3$. The
boundary $\partial M$ of $M$ consists of the points not having
neighbourhoods homeomorphic to $\matR^3$, and it is a (possibly
disconnected) surface.

\begin{prop}\label{3mfld:bd:prop}
Fix the notation of Proposition~\ref{3mfld:prop} and let
$M(\varphi)$ be the space obtained from the truncated octahedron
$O_{\rm t}$ by the gluing induced from $\varphi$. Then $M(\varphi)$
is a $3$-manifold with boundary. Moreover:
\begin{itemize}
\item There is a natural correspondence between the components of $\partial M(\varphi)$ and
the elements of $V(\varphi)$;
\item A point of $V(\varphi)$ is a $3$-manifold point of $X(\varphi)$ if and
only if the corresponding component of $\partial M(\varphi)$ is a sphere.
\end{itemize}
\end{prop}

\begin{proof}
If we remove from $O_{\rm t}$ its 6 square faces, we get a space homeomorphic
to $O\setminus O^{(0)}$. This easily implies that the space ${\rm Int}(M(\varphi))$ obtained from
$M(\varphi)$ by removing the projections of the
squares is homeomorphic to $Y(\varphi)$, so it is a 3-manifold without boundary.
The complement in $M(\varphi)$ of ${\rm Int}(M(\varphi))$ is obtained from the squares by gluing
together in pairs their edges, so it is a surface $\partial M(\varphi)$. It is
now immediate to check that ${\rm Int}(M(\varphi))$ and $\partial M(\varphi)$ fit together nicely
to give $M(\varphi)$ the structure of a 3-manifold with boundary.

By the very construction, a neighbourhood in $X(\varphi)$ of an element of $V(\varphi)$
is the cone over a component of $\partial M(\varphi)$. The conclusion easily
follows from the remark that the cone over a surface is a 3-manifold
if and only if the surface is the sphere.
\end{proof}

\section{Partial answers from lower dimension}
To face Question~\ref{basic:ques} the first issue is of course to
determine how many gluing patterns can be chosen for the faces of
the octahedron. The situation is almost but not quite similar to the
2-dimensional case of $P_k$. On one hand, to determine a gluing, we
must again start by choosing which faces get glued to which, and
this gives a wealth of $(8-1)!!=105$ choices. But on the other hand
this first information, which is sufficient for $P_k$ to determine
the gluing, is not sufficient for $O$, because there are three
distinct orientation-reversing simplicial homeomorphisms between two
given faces (obtained from each other by pre- or post-composition
with a rotation). The total number of gluing patterns is therefore $105\cdot 3^4=8505$.

Since combinatorially equivalent gluing patterns of course give rise to
homeomorphic glued spaces, the next task is to let the symmetry group of the
octahedron act on the 8505 gluings. The group has 48 elements but several
gluing patterns have themselves non-trivial symmetries, so it turns out that
the actual number (that we have determined by computer, again using
some Haskell code) is significantly larger
than $8505/48=177.1875$:

\begin{prop}\label{patterns:prop}
There exist $298$ combinatorially inequivalent patterns
of orientation-reversing gluings of the faces of $O$.
\end{prop}

Having carried out this preliminary combinatorial work, to attack
Question~\ref{basic:ques} we now need to start facing topological
issues. But the discussion of
the previous section suggests that, from the point of view
of the topology of 3-manifolds, the following is a more natural problem
than the original one:

\begin{ques}\label{bd:ques}
How many different $3$-manifolds with boundary $M(\varphi)$ can be obtained from the
truncated octahedron $O_{\rm t}$ under a gluing of its hexagonal faces
induced by some simplicial orientation-reversing gluing $\varphi$ of the
faces of the octahedron $O$?
\end{ques}

Two homeomorphic 3-manifolds of course have homeomorphic boundaries,
and we have seen above that identifying the topological type of a
surface is easy, so the first natural thing to do to attack this
question is to subdivide the inequivalent gluings $\varphi$
according to the type of $\partial M(\varphi)$.
The figures (obtained by computer) are given in
Table~\ref{bd:tab}, where $S$ is the sphere, $T$ is the torus and
$\Sigma_g$ is the surface of genus $g\geqslant 2$. (Of course we
could have denoted $S$ and $T$ by $\Sigma_0$ and $\Sigma_1$
respectively. However we already know from
Proposition~\ref{3mfld:bd:prop} that $S$ has a totally peculiar
status among surfaces from the viewpoint of topology, whereas we
will see in Section~\ref{hyp:sec} that $T$ is peculiar from the
viewpoint of hyperbolic geometry. So we prefer to keep the notation
clearly distinct).

\begin{table}
\begin{center}
\begin{tabular}{c|c}
$\partial M(\varphi)$ & $\#$(inequivalent $\varphi$'s) \\ \hline\hline
$S$ &  23\\ \hline
$S\sqcup S$ & 8\\ \hline
$S\sqcup S\sqcup S$ & 4\\ \hline
$S\sqcup S\sqcup S\sqcup S$ & 1\\ \hline
$S\sqcup S\sqcup S\sqcup S\sqcup S$ & 1\\ \hline
$T$ & 67\\ \hline
$T\sqcup S$ & 10\\ \hline
$T\sqcup S\sqcup S$ & 4\\ \hline
$T\sqcup T$ & 8 \\ \hline
$T\sqcup T\sqcup S$ & 1 \\ \hline
$\Sigma_2$ & 113 \\ \hline
$\Sigma_2\sqcup T$ &  2\\ \hline
$\Sigma_3$  & 56  \\ \hline\hline
Total & 298 \\
\end{tabular}
\end{center}
\mycap{Numbers of inequivalent gluings $\varphi$ of the faces of $O$,
subdivided according to the topological type of $\partial M(\varphi)$\label{bd:tab}}
\end{table}

Table~\ref{bd:tab} shows that the type of $\partial M(\varphi)$
gives rather limited
information towards the answer to Question~\ref{bd:ques}. In
addition this question is indeed not equivalent to our original one,
but a small variation on it turns out to be so. This will follow
from the next easy result, to state which we introduce a little
extra notation. Given a 3-manifold $M$ we will denote by $\hatM$ the
space obtained from $M$ by attaching a 3-disc to each component of
$\partial M$ homeomorphic to the sphere $S$. Note that $\hatM$ is
again a 3-manifold with boundary, but if $\partial M$ consists of
spheres only then $\partial\hatM$ is empty, \emph{i.e.}~$\hatM$ is a
closed 3-manifold.

\begin{prop}\label{X:hat:prop}
$X(\varphi)$ and $X(\varphi')$ are homeomorphic to each other if and
only if $\hatM(\varphi)$ and
$\hatM(\varphi')$ are homeomorphic to each other.
\end{prop}

In the rest of the paper we will consider the following
reformulation of Question~\ref{basic:ques}:

\begin{ques}\label{bd:cap:ques}
How many different $3$-manifolds $\hatM(\varphi)$ exist as $\varphi$
varies among the simplicial orientation-reversing gluings of the
faces of the octahedron $O$?
\end{ques}

The contribution of 2-dimensional topology towards the answer to
this question is easily deduced from Table~\ref{bd:tab} and
summarized in Table~\ref{bd:cap:tab}.

\begin{table}
\begin{center}
\begin{tabular}{c|c}
$\partial \hatM(\varphi)$ & $\#$(inequivalent $\varphi$'s) \\ \hline\hline
$\emptyset$ &  37 \\ \hline
$T$ & 81\\ \hline
$T\sqcup T$ & 9\\ \hline
$\Sigma_2$ &  113\\ \hline
$\Sigma_2\sqcup T$ &  2\\ \hline
$\Sigma_3$  & 56  \\ \hline\hline
Total &  298 \\
\end{tabular}
\end{center}
\mycap{Numbers of inequivalent gluings $\varphi$ of the faces of $O$,
subdivided according to the topological type of $\partial \hatM(\varphi)$\label{bd:cap:tab}}
\end{table}

\section{The magic of hyperbolic\\ geometry in three dimensions}\label{hyp:sec}
According to Table~\ref{bd:cap:tab}, the answer to
Question~\ref{bd:cap:ques}, and hence to Question~\ref{basic:ques},
could be any number between 6 and 298, so the information provided
by the boundary is very weak. To go further we will turn from purely
combinatorial and topological methods to geometric
ones. To this end we recall that a
Riemannian metric~\cite{Riemannian} on a manifold $M$ (possibly with
non-empty boundary) is a smoothly varying inner product on the
tangent spaces at the points of $M$. Given such a metric on $M$ one can:

\begin{itemize}
\item Define the \emph{length} $L(\alpha)$ of a smooth path $\alpha$ in $M$ by integrating the norm of
the tangent vector $\dot\alpha$;
\item Introduce on $M$ the \emph{distance} $d$ between any two points, given by the
infimum of the lengths of smooth paths joining them;
\item Define a \emph{geodesic} as a smooth path locally realizing the distance between the points
it contains;
\item Define a submanifold $N$ of $M$ (for instance, the boundary $\partial M$) to
be \emph{totally geodesic} if any geodesic meeting $N$
in more than one point is actually contained in $N$; for the sake of brevity
we will henceforce omit the word ``totally'';
\item Define the (possibly infinite) volume of $M$ by integrating the
square root of the determinant of the matrix representing the metric
in local charts;
\item Assign to each 2-plane $P$ contained in a tangent space $T_xM$ a real
number $K_x(P)$ called the \emph{sectional curvature} of $M$ at $x$ along $P$.
\end{itemize}

We recall that for a Riemannian 2-manifold, \emph{i.e.}~a surface,
the curvature at a point $x$ measures the extent to which the
circle of a very small radius $r$ centred at $x$ fails to have
length $2\pi r$, according to the formula
$$K_x=\lim\limits_{r\to 0^+}6\cdot\frac{2\pi r-L(C_x(r))}{r^2},\qquad C_x(r)=\{y:\ d(x,y)=r\}.$$
(The 6 is introduced in the limit so that the round sphere of radius 1 in
$\matR^3$ with metric induced from $\matR^3$ has curvature $+1$).
The sectional curvature $K_x(P)$ is the curvature at $x$ of a germ of geodesic
surface with tangent plane $P$ at $x$.

We will now restrict to the type of metrics that turns out
to be mostly relevant for 3-dimensional topology.
To this end we consider an orientable manifold $M$ with compact boundary $\partial M$
(but possibly non-compact itself). We define a Riemannian
metric on $M$ to be \emph{hyperbolic} if:

\begin{itemize}
\item The sectional curvature is $-1$ at every point and along every tangent
plane;
\item $M$ is a complete metric space with respect to the distance
induced by the Riemannian metric;
\item The volume of $M$ is finite.
\end{itemize}

A number of spectacular results have been proved about hyperbolic manifolds in dimension 3
over the last 30 years, many of which due to or inspired by Bill Thurston~\cite{bible}.
We list here those we will need to refer to, but we mention that there are very many more:

\begin{itemize}

\item \textbf{Mostow's rigidity theorem}~\cite{bible,ratcliffe,lectures}:
\emph{Any two hyperbolic metrics on the same manifold are isometric
to each other}. This implies in particular that any hyperbolic
invariant, such as the volume or the length of the shortest
geodesic, is automatically a topological invariant, which gives an
extremely powerful method for distinguishing
manifolds.
Note that the assumption that the volume should be finite is crucial
for rigidity;

\item \textbf{Cusps and compactification}~\cite{bible,ratcliffe,lectures}:
\emph{If $M$ is hyperbolic then $M$ decomposes as $M_0\cup C_1\cup\ldots\cup C_k$, where $M_0$ is
a compact manifold with boundary $\partial M_0$ consisting of $\partial M$ together with
$k$ tori $T_1,\ldots,T_k$, and $C_i=T_i\times [0,\infty)$ with $T_i=T_i\times\{0\}$.
In particular, $M$ compactifies to a manifold $\Mbar$ by adding the tori
$T_i\times\{\infty\}$.} Since each boundary component of $M$ has a metric of constant curvature $-1$,
the Gauss-Bonnet formula~\cite{Riemannian} implies that it cannot be the sphere
or the torus. So $\Mbar$ is a compact manifold without boundary spheres and $M$
is obtained from $\Mbar$ by removing all the boundary tori;

\item \textbf{Thurston's hyperbolization proved by
Perelman}~\cite{bible,perelman,perelman2,perelman3}:
\emph{If $\Mbar$ is a compact orientable manifold without boundary
spheres and different from the solid torus, $\pi_1(M)$ is an
infinite group, and $M$ does not contain any essential embedded
sphere, disc, annulus or immersed torus, then the manifold
$M=\Mbar\setminus$(boundary tori) is hyperbolic}. The notion of
``essentiality'' of an embedded surface is too technical to be
reproduced here, but it basically means that cutting $\Mbar$ along
the surface one gets some non-obvious simplification of $\Mbar$; an
immersed torus is essential if its fundamental group injects in that
of $\Mbar$ but none of its conjugates is contained in a subgroup of
$\pi_1(M)$ arising from a toric boundary component. Combining two
well-known general theorems (the Haken-Kneser-Milnor decomposition
along spheres~\cite{Hempel} and the Jaco-Shalen-Johansson
decomposition along tori~\cite{matbook}) one knows that a set of
``building blocks'' for all 3-manifolds is given by:
\begin{itemize}
\item those as described in the statement of the hyperbolization theorem;
\item those finitely covered by a closed
simply connected manifold;
\item those of a type that has been classified long ago (Seifert manifolds).
\end{itemize}
So Thurston's hyperbolization theorem together with
the positive solution to the Poincar\'e conjecture (also proved by Perelman)
can be interpreted as a classification result for all 3-manifolds, and shows
that hyperbolic geometry plays a central r\^ole in 3-dimensional topology;

\item \textbf{Epstein-Penner canonical decomposition}~\cite{EP}:
\emph{If $M$ is hyperbolic, $\partial M=\emptyset$ and $M$ is non-compact then
$M$ has a canonical realization as a gluing of ideal polyhedra}. An ideal polyhedron
is just one with its vertices removed, and we have already met
in Proposition~\ref{3mfld:prop} the idea of gluing together such objects.
The statement means that the geometry (and hence the topology,
thanks to rigidity) of $M$ determines not only the number of polyhedra and
their type but also the combinatorics of the gluing pattern;

\item \textbf{Kojima canonical decomposition}~\cite{koji1,koji2}:
\emph{If $M$ is hyperbolic and $\partial M\neq\emptyset$ then
$M$ has a canonical realization as a gluing of truncated/ideal polyhedra}.
This statement has precisely the same meaning as the previous one,
except that in this case some vertices of the polyhedra are removed
and give rise to the cusps of $M$, while other vertices are truncated,
and the truncation polygons glue together to give the geodesic boundary of $M$;

\item \textbf{Algorithmic approach}:

\begin{itemize}

\item \textbf{Search for hyperbolic structure}~\cite{bible,SnapPea,lectures,akira,fp,Orb}:
\emph{Given a manifold $\Mbar$ realized as a gluing of truncated
tetrahedra, one can try to build the hyperbolic structure on
$M=\Mbar\setminus$(boundary tori) by choosing for the tetrahedra
geometric shapes (parameterized by numbers) and imposing these
structures to glue together nicely (which translates into equations
in the parameters)}. This searching method for the hyperbolic
structure is not guaranteed to always work, but in practice it most
often does if one starts from a realization having the minimal
possible number of tetrahedra. The algorithm was implemented in the
case of empty boundary by Weeks and in the case of non-empty
boundary by Frigerio-Martelli-Petronio and by Heard;

\item \textbf{Search for Epstein-Penner decomposition}~\cite{SaWe,SnapPea}:
\emph{The can\-onical realization of a cusped hyperbolic manifold with
empty boun\-dary can be searched for
algorithmically starting from a geometric realization as a gluing of
ideal tetrahedra}. This method is based on a certain ``tilt formula'' due
to Sakuma and Weeks. Again the algorithm is not guaranteed to always
converge but in practice it does. It was implemented by Weeks;

\item \textbf{Search for Kojima decomposition}~\cite{akira,fp,Orb}:
\emph{The canonical realization of a cusped hyperbolic manifold with
non-empty boundary can be searched for algorithmically starting from a
geometric realization as a gluing of truncated/ideal tetrahedra}. This method is
based on a variation of the ``tilt formula'' established by
Ushijima. Once more the algorithm is not guaranteed to always
converge but in practice it does; it was implemented by
Frigerio-Martelli-Petronio and by Heard.

\end{itemize}
\end{itemize}

\section{Answers from hyperbolic geometry}\label{ans:sec}
In this section we will employ the rich technology of hyperbolic geometry
to show that the combinatorially
inequivalent gluing patterns of Table~\ref{bd:cap:tab}
actually give rise to very many different manifolds. We begin
with the following striking fact, that was initially discovered
from a computer experiment~\cite{fmp3} and was later proved theoretically:

\begin{prop}\label{56:prop}
The $56$ gluing patterns of Table~\ref{bd:cap:tab} giving rise to a
boundary of genus $3$ define $56$ pairwise distinct hyperbolic
manifolds with geodesic boundary. These manifolds all have the same
volume $11.448776...$.  The homology is $\matZ^3$ for $52$ of them
and $\matZ_3\times\matZ^3$ for four of them.
\end{prop}

\begin{proof}
An easy computation of Euler characteristic shows that a gluing $\varphi$ defines
a manifold $\hatM(\varphi)$ bounded by $\Sigma_3$ if and only if it identifies
all 12 edges to each other. We want to show that such an $\hatM(\varphi)$ is hyperbolic
with geodesic boundary by choosing a hyperbolic shape of the truncated
octahedron $O_{\rm t}$ that is matched by $\varphi$. Since all edges
are glued together, this can only happen if the geometric shape is such that
all edges have the same length, \emph{i.e.}~$O_{\rm t}$ is regular.
If this is the case, all dihedral angles are also the same, so they must
all be $2\pi/12$. Such an octahedron certainly does not exist in Euclidean
or spherical geometry, but it does in hyperbolic geometry. This implies
that $\hatM(\varphi)$ is indeed hyperbolic.

Let us now analyze the Kojima canonical decomposition of $\hatM(\varphi)$.
To this end we recall~\cite{koji1,koji2} that it is dual
to the cut locus of the boundary, \emph{i.e.}~to the set of points having
multiple shortest paths to $\partial \hatM(\varphi)$. Using the
fact that $\hatM(\varphi)$ is the gluing of a regular truncated octahedron, which
is totally symmetric, it is not too difficult to show that the Kojima
decomposition is given by the octahedron itself, with its gluing
pattern $\varphi$. This implies that the geometry of $\hatM(\varphi)$, and hence its
topology, determines $\varphi$. Therefore different $\varphi$'s give
rise to different $\hatM(\varphi)$'s.

The fact that all the $\hatM(\varphi)$'s have the same volume is obvious.
For the computation of the value of this volume and for the
determination of the homology we address the reader to~\cite{fmp3}.
\end{proof}

Turning to the other gluing patterns in Table~\ref{bd:cap:tab}, we
now cite some facts that follow from the
results in~\cite{fmp3}:

\begin{prop}\label{Sigma2:prop}
The two $\varphi$'s such that $\partial \hatM(\varphi)=\Sigma_2\sqcup T$
define distinct hyperbolic manifolds with geodesic boundary and one cusp.
The group of $113$ inequivalent $\varphi$'s such that $\partial \hatM(\varphi)=\Sigma_2$
contains:
\begin{itemize}
\item A group of $14$ $\varphi$'s defining $14$ distinct hyperbolic
manifolds with geodesic boundary and Kojima decomposition given
by a truncated octahedron;

\item A group of $8$ $\varphi$'s defining $8$ distinct hyperbolic
manifolds with geodesic boundary and Kojima decomposition given
again by a truncated octahedron, but with a
different geometric shape;

\item A group of $4$ $\varphi$'s defining $4$ distinct hyperbolic
manifolds with geodesic boundary and Kojima decomposition given
by two square pyramids.
\end{itemize}
The manifolds arising from these three groups of gluings are
all distinct from each other. Moreover any other
hyperbolic $\hatM(\varphi)$ with $\partial \hatM(\varphi)=\Sigma_2$
has Kojima decomposition consisting of tetrahedra only.
\end{prop}

The gluing patterns in Table~\ref{bd:cap:tab} not covered by the
results just stated have been analyzed by Damian Heard using his
software ``Orb'', with the following results:

\begin{prop}\label{damian:prop}
\begin{itemize}
\item Consider the set of $87$ inequivalent $\varphi$'s such that $\partial \hatM(\varphi)=\Sigma_2$
but not included in Proposition~\ref{Sigma2:prop}. Then this set contains at least
$37$ elements defining hyperbolic manifolds with geodesic boundary,
and these manifolds are all distinct from each other;
\item The set of $81$ inequivalent $\varphi$'s with $\partial \hatM(\varphi)=T$
contains at least $11$ elements defining one-cusped hyperbolic
manifolds, and there are $9$ different manifolds arising from these
$11$ $\varphi$'s;
\item The set of $9$ inequivalent $\varphi$'s with $\partial \hatM(\varphi)=T\sqcup T$
contains at least $2$ elements defining two-cusped hyperbolic
manifolds, and these two manifolds are distinct.
\end{itemize}
\end{prop}

As one sees, hyperbolic geometry only gives a partial answer to
Question~\ref{bd:cap:ques} when the boundary
is not $\Sigma_3$. In addition the octahedron is often not related
to the geometry of the result. For this reason we will consider a
more general question in the next section, but the appendix contains
the full answer to the original question on the gluings of the
octahedron.

\section{Other enumeration problems and results}\label{enu:sec}
Since the octahedron can be subdivided (in three different ways) into four tetrahedra,
listing the different manifolds arising from gluings of the octahedron is a special
case of the following:

\begin{ques}\label{c:pre:ques}
Given $c>0$, how many different $3$-manifolds can be obtained by gluing together in pairs
the faces of a disjoint union of $c$ tetrahedra?
\end{ques}

As usual we consider this question only under the restriction that the face-pairings
should be orientation-reversing, and we remove from the space obtained
after the gluing an open regular neighbourhood of the non-manifold points, so to get
compact manifolds bounded by surfaces of positive genus.

A systematic approach to Question~\ref{c:pre:ques} was proposed by
Matveev~\cite{Acta}, who introduced a certain \emph{complexity
theory}. To describe it, we need to recall the following
definitions and facts~\cite{Hempel}:
\begin{itemize}
\item A $3$-manifold is called \emph{prime} if every embedded 2-sphere it contains is the boundary
of an embedded 3-disc;
\item The connected sum $\#$ of two $3$-manifolds is obtained by removing an open 3-disc
from each and gluing together the boundary spheres thus created;
\item Every $3$-manifold can be expressed in a unique way as a connected sum of a finite
number of prime 3-manifolds.
\end{itemize}
With this in mind, Matveev noticed that a good notion of complexity for
a 3-manifold $M$, \emph{i.e.}~a measure $c(M)\in\matN$ of how complicated $M$
is, should have the property of being additive under connected sum.
This is however certainly not the case if one defines $c(M)$
as the minimal number of tetrahedra needed to realize $M$ via face-pairings,
because $S^3$ would have positive complexity, but $c(M\#S^3)=M$ for any $M$.
This problem was solved in~\cite{Acta} by defining $c(M)$ through a different
and more flexible object than a decomposition into tetrahedra, called a \emph{simple spine}.
We will not recall what such a thing is in this paper, but we will mention
that among simple spines there are some having additional properties,
called \emph{special spines}, that turn out to be perfectly equivalent
to decompositions into tetrahedra. The main features of complexity theory,
that we state for the case of closed orientable manifolds but have extensions to
more general situations, are now the following:
\begin{itemize}
\item $c(M)$ is a non-negative integer;
\item $c(M_1\#M_2)=c(M_1)+c(M_2)$;
\item The prime 3-manifolds $M$ with $c(M)=0$ are the sphere $S^3$, the projective space $\matR\matP^3$,
the lens space $L(3,1)$, and the product $S^2\times S^1$;
\item If $M$ is prime and $c(M)>0$ then $c(M)$ is precisely the minimal number of
tetrahedra needed to realize $M$ via face-pairing.
\end{itemize}

>From the viewpoint of complexity theory, the following
variation on Question~\ref{c:pre:ques} is therefore more
natural:

\begin{ques}\label{c:ques}
Given $c>0$, how many different prime
$3$-manifolds can be obtained by gluing together in pairs
the faces of a disjoint union of $c$ tetrahedra?
\end{ques}

This question was investigated both theoretically and using computers by Matveev, Martelli-Petronio,
Martelli and Matveev (see~\cite{bruno:surv} for a list of the relevant papers and websites).
The answer was obtained for $c\leqslant 12$ and it is described in Table~\ref{closed:tab},
where the information of how many hyperbolic manifolds were found is also added.
We further note that the first four such manifolds, appearing in complexity 9 and
first discovered in~\cite{mp1},
are precisely those with the four smallest known volumes. This supports
a conjecture of Matveev and Fomenko that low volume should appear in low complexity.

\begin{table}
\begin{center}
\begin{tabular}{l||c|c|c|c|c|c|c|c|c|c|c|c|c}
$c$ & 0 & 1 & 2 & 3 & 4 & 5 & 6 & 7 & 8 & 9 & 10 & 11 & 12 \\
\hline\hline tot & 4 & 2 & 4 & 7 & 14 & 31 & 74 & 175 & 436 &
1154 & 3078 & 8343 & 23431\\ \hline
hyp & -- & -- & -- & -- & -- & -- & -- & -- & -- & 4 & 25 & 120 & 459 \\
\end{tabular}
\end{center}
\mycap{Numbers of prime closed manifolds up to complexity 12\label{closed:tab}}
\end{table}

We have spoken so far of closed manifolds, but a variation of complexity theory applies
to bounded compact 3-manifolds too. However in this case the equality between $c(M)$
and the minimal number of tetrahedra in a decomposition requires $M$ not only to
be prime but also to have incompressible boundary. We will not review this notion
here, but we mention it is implied by hyperbolicity. Using this fact, two groups of
people investigated hyperbolic 3-manifolds of low complexity.
Callahan-Hildebrandt-Weeks~\cite{CaHiWe}, using SnapPea~\cite{SnapPea},
classified cusped manifolds up to complexity 7,
getting the figures of Table~\ref{cusped:tab}.
(A census exists also for non-orientable cusped hyperbolic manifolds, and
it contains much fewer elements).

\begin{table}
\begin{center}
\begin{tabular}{r||c|c|c|c|c|c|c}
complexity & 1 & 2 & 3 & 4 & 5 & 6 & 7 \\ \hline\hline
one cusp & -- & 2 & 9 & 52 & 223 & 913 & 3388 \\ \hline
two cusps & -- & -- & -- & 4 & 11 & 48 & 162 \\ \hline
three cusps & -- & -- & -- & -- & -- & 1 & 2 \\
\end{tabular}
\end{center}
\mycap{Numbers of hyperbolic cusped manifolds up to complexity 7\label{cusped:tab}}
\end{table}

Frigerio-Martelli-Petronio~\cite{fmp3} (see also ~\cite{fmp1,fmp2})
more recently turned their attention to
hyperbolic manifolds with non-empty geodesic boundary, and found the numbers
shown in Table~\ref{geo:bd:tab}.

\begin{table}
\begin{center}
\begin{tabular}{r||c|c|c}
complexity & 2 & 3 & 4 \\ \hline\hline
$\Sigma_2$ & 8 & 76 & 628 \\ \hline
$\Sigma_2\cup T$  & -- & 1 & 18 \\ \hline
$\Sigma_2\cup T\cup T$  & -- & -- & 1 \\ \hline
$\Sigma_3$  & -- & 74 & 2034 \\ \hline
$\Sigma_3\cup T$  & -- & -- & 12 \\ \hline
$\Sigma_4$  & -- & -- & 2340 \\
\end{tabular}
\end{center}
\mycap{Numbers of hyperbolic manifolds with compact and non-empty geodesic
boundary up to complexity 4, subdivided according to the type of $\partial\Mbar$\label{geo:bd:tab}}
\end{table}

\section*{Appendix (by Damian Heard and\\ Ekaterina Pervova): The complete answer}
The complete answer to Question~\ref{basic:ques} (actually, to its
equivalent formulation Question~\ref{bd:cap:ques}) is contained in
Table~\ref{final:tab}. The hyperbolic entries of this table have
been discussed in the body of the paper (note that $63$ arises as
the sum of $14+8+4$ from Proposition~\ref{Sigma2:prop} plus the $37$
from Proposition~\ref{damian:prop}). They were all obtained using
the software ``Orb''~\cite{Orb} (or confirmed by ``Orb'' when they
were known from~\cite{fmp3}). Using randomization of triangulations,
``Orb'' also grouped together homeomorphic manifolds it could not
construct a hyperbolic structure for, finding precisely the numbers
shown in the table. To show that the numbers were indeed the correct
ones we were then left with the following tasks:

\begin{enumerate}
\item Prove that among the groups of apparently non-hyperbolic manifolds
produced by ``Orb'' there were no duplicates;
\item Prove that the apparently non-hyperbolic manifolds indeed
were non-hyperbolic.
\end{enumerate}

\begin{table}
\begin{center}
\begin{tabular}{c||c|c|c||c}
boundary type & $\#$(gluings) & hyperbolic & non-hyperbolic & total
\\ \hline\hline $\emptyset$ &  37 & -- & 17 & 17\\ \hline $T$ &  81
& $9$ & $21$ & 30\\ \hline $T\sqcup T$ & 9 & 2 & 5 & 7\\ \hline
$\Sigma_2$ &  113 & 63 & 16 & 79\\ \hline $\Sigma_2\sqcup T$ & 2 &
2 & -- & 2 \\ \hline $\Sigma_3$  & 56 & 56  & -- & 56 \\
\hline\hline
Total & 298 &  132 &  59 &  191  \\
\end{tabular}
\end{center}
\mycap{Numbers of distinct manifolds arising
from orientation-preserving gluings of the faces of a truncated octahedron,
with boundary spheres capped off\label{final:tab}}
\end{table}

\paragraph{The ``Recognizer''}
The basic tool we have used to achieve both these tasks is the
``3-Manifold Recognizer,'' a software written by Tarkaev and
Matveev~\cite{recognizer}. The input to this program is a
triangulation of a 3-manifold $M$ (so we had to subdivide our
octahedron, and actually switch to the dual viewpoint of spines) and
its output is the ``name'' of $M$, by which we mean the following:
\begin{itemize}
\item For a Seifert $M$, (one of) its Seifert structure(s);
\item For a hyperbolic $M$, its presentation(s) as a Dehn filling of a manifold in the
tables of Weeks~\cite{CaHiWe};
\item For a prime $M$ having JSJ decomposition into more than block,
the name (as just illustrated) of the blocks, together with the
gluing instructions between the blocks;
\item For a non-prime manifold, the names (as just illustrated) of its
prime summands.
\end{itemize}
The program is not guaranteed to always find the name of the
manifold (for instance, it does not even attempt to do this for
manifolds with boundary
of genus 2 or more, and it happens to fail also in other cases).
But it can always compute the first homology
and, in the case of boundary of genus at most $1$,
the Turaev-Viro invariants~\cite{TV}, which turned out to be very useful for us.

\paragraph{Closed manifolds}
For the 37 gluing patterns in Table~\ref{bd:cap:tab} defining closed
manifolds, task 2 (proof of non-hyperbolicity) was not an issue,
since it has been known for a long time~\cite{matbook} that closed
hyperbolic manifolds start appearing in complexity 9, whereas a
gluing of the octahedron has complexity at most 4. For task 1 (proof
that the grouping by ``Orb'' contains no duplicates), we have run
the ``Recognizer,'' that successfully identified all the manifolds
(this was also independently done by Tarkaev). From the names (all
manifolds turned out to be Seifert or connected sums of Seifert) we
could see there were either 16 or 17 of different ones. The only
uncertainty was related to the fact that two manifolds were
recognized to be the connected sum of two copies of the lens space
$L(3,1)$, the point being that $L(3,1)$ has no orientation-reversing
automorphism. So, even if one looks at orientable but non-oriented
manifolds, there are two distinct ways of performing the connected
sum of $L(3,1)$ with itself.

We then had to examine the two triangulations by hand, introducing
an arbitrary orientation on each of them. For each triangulation we
then found the non-trivial sphere realizing the connected sum.
Cutting along this sphere and capping off, we saw that in one case
the two connected summands were distinctly oriented copies of
$L(3,1)$, while in the other case we got two copies of $L(3,1)$ with
the same orientation. Thus, the manifolds were respectively
$L(3,1)\#(-L(3,1))$ and $L(3,1)\# L(3,1)$, hence distinct.

\paragraph{One-cusped manifolds}
For the 81 triangulations from Table~\ref{bd:cap:tab} giving
manifolds bounded by one torus we followed approximately the same
strategy.  After the work with ``Orb''
already described above (recognition of hyperbolic manifolds
and grouping of the others), we had 21 manifolds
$\{M_i\}_{i=1}^{21}$ we should realize tasks 1 and 2 on. We
proceeded as follows:
\begin{enumerate}
\item To carry out task~1, \emph{i.e.} to prove that there were no
duplicates among the 21 manifolds, we again employed the
``Recognizer'', using which we calculated their first homology
groups and Turaev-Viro invariants up to order~16. From this
computation we deduced that $M_i\not\cong M_j$ for $1\leqslant
i<j\leqslant 21$ except possibly for $i=1,2,3,4$ and $j=i+4$. For
the four pairs of manifolds left, we showed the homeomorphism was
impossible by analyzing the JSJ decompositions. Specifically, $M_1$
and $M_5$ turned out to be Seifert and distinct, and the same
happened for $M_2$ and $M_6$, whereas $M_3$ and $M_7$ had
non-trivial JSJ decompositions, with the same blocks but different
gluing matrices, and analogously for $M_4$ and $M_8$.

\item In order to complete task~2, \emph{i.e.} to show that the
manifolds were not hyperbolic, we computed the JSJ decomposition also for
$M_9,\ldots,M_{21}$ (again using the ``Recognizer''). In all
cases we obtained a decomposition (sometimes trivial) consisting
of Seifert pieces. (In one case the ``Recognizer'' failed to return
the answer right away, but we were able to transform the
triangulation by hand into one that the ``Recognizer'' could
handle).
\end{enumerate}

\paragraph{Two-cusped manifolds}
``Orb'' reduced the issue of counting these manifolds to that of
realizing tasks 1 and 2 on 5 of them, which was easy in this case
using the ``Recognizer'':

\begin{enumerate}
\item For task~1, we
calculated the Turaev-Viro invariants of the manifolds, which was
sufficient to prove that all 5 manifolds were indeed
distinct;

\item For task~2, we found the
JSJ decomposition of the manifolds. It turned out that either this
decomposition was non-trivial or the manifold in question was
Seifert. Hence indeed none of the manifolds was hyperbolic.
\end{enumerate}

\paragraph{Genus-2 boundary}
For the 113 gluing patterns from Table~\ref{bd:cap:tab} giving manifolds
bounded by one genus-2 surface, ``Orb'' showed there were 16
apparently non-hyperbolic and apparently distinct manifolds, and
also provided a presentation for the fundamental group of
each of them. These 16 manifolds were then treated as follows:

\begin{enumerate}
\item To address task~1, we computed the first homology
of the manifolds (with the ``Recognizer'') and the Turaev-Viro
invariants (by hand). Using these invariants, we were able to break
down the set of 16 manifolds into some pairs and some 4-tuples of
potentially equal manifolds, and also to single out two individual
manifolds which were distinct from all the other ones. For the
resulting groups of manifolds we then used the presentations of
their fundamental groups to calculate the homology of the 3-fold
coverings, which allowed us to distinguish almost all manifolds in
these groups. This left only two pairs of manifolds, to which the
same technique could not apply because the fundamental groups were
isomorphic. To tackle these, we found, by hand, their JSJ
decompositions and examined the gluings. In both cases the
decomposition consisted of a genus-2 handlebody and a solid torus
glued along an annulus. In one case the annuli used in the gluing
were the same but the two gluings were along non-isotopic
homeomorphisms. In the second case the annuli in the boundaries of
the respective handlebodies were different, \emph{i.e.} not related
by any homeomorphism of the handlebodies.

\item To realize task~2, we showed that each of the 16 manifolds
contained an essential annulus (which was explicitly constructed).
For this, we employed the dual viewpoint of so-called \emph{simple
spines}, constructing for each of the manifolds in question a simple
spine with an annulus component or a M\"obius strip component. We
then checked that in each case the annulus dual to the core of the
annulus component or that dual to the boundary of the M\"obius strip
component, whichever was relevant, was actually essential in the
manifold.
\end{enumerate}

\vspace{1cm}

\noindent Dipartimento di Matematica Applicata\\
Via Filippo Buonarroti, 1C\\
56127 PISA -- Italy\\
{\tt petronio@dm.unipi.it}

\vspace{.5cm}

\noindent RedTribe\\
Carlton, Victoria\\
Australia 3053\\
{\tt damianh@redtribe.com}

\vspace{.5cm}

\noindent Dipartimento di Matematica Applicata\\
Via Filippo Buonarroti, 1C\\
56127 PISA -- Italy\\
{\tt pervova@csu.ru}


\begin{thebibliography}{99}

\bibitem{Riemannian}
\textsc{W.~Klingenberg},
``Riemannian Geometry,''
De Gruyter Std. in Math. Vol. 1, Berlin-New York, 1982.

\bibitem{lectures}
\textsc{R.~Benedetti -- C.~Petronio}, ``Lectures on Hyperbolic
Geometry,'' Springer-Verlag, Berlin-Heidelberg-New York, 1992.

\bibitem{CaHiWe}
\textsc{P.~J.~Callahan -- M.~V.~Hildebrandt -- J.~R.~Weeks}, \emph{A
census of cusped hyperbolic 3-manifolds. With microfiche
supplement}, Math. Comp. \textbf{68} (1999), 321-332.

\bibitem{EP}
\textsc{D.~B.~A.~Epstein, R.~C.~Penner}, \emph{Euclidean
decomposition of non-compact hyperbolic manifolds}, J. Differential
Geom. (1) \textbf{27} (1988), 67-80.


\bibitem{fmp1}
\textsc{R.~Frigerio -- B.~Martelli -- C.~Petronio},
\emph{Complexity and Heegaard genus of an infinite class of
compact $3$-manifolds}, Pacific J. Math. \textbf{210} (2003),
283-297.

\bibitem{fmp2}
\textsc{R.~Frigerio -- B.~Martelli -- C.~Petronio}, \emph{Dehn
filling of cusped hyperbolic $3$-manifolds with geodesic
boundary}, J. Differential Geom. \textbf{64} (2003), 425-455.

\bibitem{fmp3}
\textsc{R.~Frigerio -- B.~Martelli -- C.~Petronio}, \emph{Small
hyperbolic $3$-manifolds with geodesic boundary}, Exp. Math.
\textbf{13} (2004), 171-184.

\bibitem{fp}
\textsc{R.~Frigerio -- C.~Petronio}, \emph{Construction and
recognition of hyperbolic $3$-manifolds with geodesic boundary},
Trans. Amer. Math. Soc. \textbf{356} (2004), 3243-3282.

\bibitem{Orb}
\textsc{D.~Heard}, \emph{``Orb''}, The computer program for finding
hyperbolic structures on hyperbolic 3-orbifolds and 3-manifolds,
available from {\tt
http://www.ms.unimelb.edu.au/$\sim$snap/orb.html}

\bibitem{Hempel}
\textsc{J.~Hempel}, ``3-Manifolds,'' Ann. of Math. Studies,
\textbf{86}, Princeton, (1976).

\bibitem{koji1}
\textsc{S.~Kojima}, \emph{Polyhedral decomposition of hyperbolic
manifolds with boundary}, Proc. Work. Pure Math. \textbf{10} (1990),
37-57.

\bibitem{koji2}
\textsc{S.~Kojima}, \emph{Polyhedral decomposition of hyperbolic
$3$-manifolds with totally geodesic boundary}, In: ``Aspects of
low-dimensional manifolds, Kinokuniya, Tokyo'', Adv. Stud. Pure
Math. \textbf{20} (1992), 93-112.

\bibitem{bruno:surv}
\textsc{B.~Martelli}, \emph{Complexity of $3$-manifolds}, In: ``Spaces
of Kleinian groups'', London Math. Soc. Lec. Notes Ser. \textbf{329}
(2006), 91-120.

\bibitem{mp1}
\textsc{B.~Martelli, C.~Petronio}, \emph{$3$-manifolds having
complexity at most $9$}, Exp. Math. \textbf{10} (2001), 207-237.

\bibitem{mp2}
\textsc{B.~Martelli, C.~Petronio}, \emph{A new decomposition theorem
for $3$-manifolds}, Illinois J. Math. \textbf{46} (2002), 755-780.

\bibitem{DAN}
\textsc{S.~V.~Matveev}, \emph{The complexity of three-dimensional
manifolds and their enumeration in the order of increasing
complexity}, Soviet Math. Dokl. \textbf{38} (1989), 75-78.

\bibitem{Acta}
\textsc{S.~V.~Matveev},  \emph{Complexity  theory  of
three-dimensional manifolds}, Acta Appl. Math. \textbf{19}
(1990), 101-130.

\bibitem{matbook}
\textsc{S.~V.~Matveev}, \emph{Algoritghmic topology and
classification of $3$-manifolds}, ACM-monographs Vol. 9,
Springer-Verlag, Berlin-Heidelberg-New York, 2003.

\bibitem{DAN05}
\textsc{S.~V.~Matveev}, \emph{Recognition and tabulation of
$3$-manifolds}, Dokl. Akad. Nauk \textbf{400} (2005), 26-28
(Russian).

\bibitem{recognizer}
\textsc{S.~V.~Matveev -- V.~V.~Tarkaev}, \emph{``Three-manifold
Recognizer''}, a computer program for recognition of 3-manifolds,
available from {\tt http://www.csu.ac.ru/$\sim$trk/spine/}.

\bibitem{perelman}
\textsc{G.~Perelman}, \emph{The entropy formula for the Ricci flow
and its geometric applications}, preprint {\tt math.DG/0211159}.

\bibitem{perelman2}
\textsc{G.~Perelman}, \emph{Ricci flow with surgery on
three-manifolds}, preprint {\tt math.DG/0303109}.

\bibitem{perelman3}
\textsc{G.~Perelman}, \emph{Finite extinction time for the solutions
to the Ricci flow on certain three-manifolds}, preprint
{\tt math.DG/0307245}.

\bibitem{ratcliffe}
\textsc{J.~G.~Ratcliffe}, ``Foundations of Hyperbolic Manifolds,''
Second Edition, Graduate Texts in Math. \textbf{149},
Springer-Verlag, New York, 2006.

\bibitem{SaWe}
\textsc{M.~Sakuma -- J.~R.~Weeks}, \emph{The generalized tilt
formula}, Geom. Dedicata \textbf{50} (1995), 1-9.


\bibitem{marty}
\textsc{M.~Scharlemann}, \emph{Heegaard splittings of compact
$3$-manifolds}, In: ``Handbook of geometric topology'' (R.~Daverman and R.~Sherr, eds.),
North-Holland, Amsterdam (2002), 921-953.


\bibitem{bible}
\textsc{W.~P.~Thurston}, ``The Geometry and Topology
of $\,3$-manifolds,'' mimeographed notes, Princeton, 1979.


\bibitem{TV}
\textsc{V.~Turaev -- O.~Viro}, \emph{State sum invariants of
$3$-manifolds and quantum 6j-symbols}, Topology (4) \textbf{31}
(1992), 865-902.

\bibitem{akira}
\textsc{A.~Ushijima}, \emph{The tilt formula for generalized
simplices in hyperbolic space}, Discrete Comput. Geom. \textbf{28}
(2002), 19-27.


\bibitem{SnapPea}
\textsc{J.~R.~Weeks}, \emph{``SnapPea''}, The hyperbolic structures
computer program, available from {\tt www.geometrygames.org}.



\end{thebibliography}
\end{document}